\documentclass[final]{amsart}
\usepackage{amssymb,amsfonts}
\usepackage{euscript,mathrsfs}
\usepackage{latexsym}
\usepackage{amscd}
\usepackage{amsmath}
\usepackage{epsf}
\usepackage[arrow,matrix]{xy} 
\usepackage{color}
\usepackage{graphicx}
\usepackage{amsthm,amsopn}
\usepackage{fixltx2e,mparhack} 
\usepackage{microtype}
\usepackage[utf8]{inputenc}
\usepackage{enumerate}
\usepackage{bbm}
\usepackage[USenglish]{babel}

\providecommand{\abs}[1]{\left\lvert#1\right\rvert} 

\newcommand{\e}{\mathrm{e}} 
\newcommand{\Ind}[1]{\mathbbm{1}_{\lbrace #1 \rbrace}} 
\newcommand{\comment}[1]{ }

\newcommand{\ie}{i.e.\;}  

\renewcommand{\subset}{\subseteq}  
\renewcommand{\supset}{\supseteq}
\newcommand{\defas}{\mathrel{\mathop{:}}=}   
\DeclareMathOperator*{\supp}{supp}  

\DeclareMathOperator*{\conv}{conv} 
\newcommand{\set}[1]{\left\lbrace #1 \right\rbrace} 
\newcommand{\fkn}{\mathcal{F}_{k,N}}

\newcommand{\akn}{\mathcal{A}_{k,N}}

\theoremstyle{plain}
\newtheorem{thm}{Theorem}

\newtheorem{prop}[thm]{Proposition}
\newtheorem*{cor}{Corollary}
\newtheorem{alg}[thm]{Algorithm}

\theoremstyle{definition}
\newtheorem{defn}[thm]{Definition}

\newtheorem{exmp}[thm]{Example}

\theoremstyle{remark}
\newtheorem*{rem}{Remark}

\begin{document}
\title{Hierarchical Models, Marginal Polytopes, and Linear Codes}

\author{Thomas Kahle}
\address{Max Planck Institute for Mathematics in the Science \\
         Inselstrasse 22 \\
         D-04103 Leipzig, Germany}
\email{kahle@mis.mpg.de}
\thanks{The first author is supported by the Volkswagen Foundation} 
\author{Walter Wenzel}
\address{Fachhochschule Magdeburg \\
         Fachbereich Wasser- und Kreislaufwirtschaft \\
         Breitscheidstrasse 2 \\
         D-39114 Magdeburg, Germany}
\email{walter@math.uni-bielefeld.de}
\author{Nihat Ay}
\address{Max Planck Institute for Mathematics in the Science \\
  Inselstrasse 22 \\
  D-04103 Leipzig, Germany \\
  and 
  Santa Fe Institute\\
  1399 Hyde Park Road \\ 
  Santa Fe, NM 87501, USA
}
\email{nay@mis.mpg.de}
\thanks{The third author is supported by the Santa Fe Institute}

\begin{abstract}
  In this paper, we explore a connection between binary hierarchical
  models, their marginal polytopes and codeword polytopes, the convex
  hulls of linear codes.\comment{Via the sufficient statistics, each
    hierarchical model is mapped to a convex polytope, the marginal
    polytope. We realize the marginal polytopes as 0/1-polytopes and
    show that their vertices form a linear code.}The class of linear
  codes that are realizable by hierarchical models is determined. We
  classify all full dimensional polytopes with the property that their
  vertices form a linear code and give an algorithm that determines
  them.
\end{abstract}

\date{\today}

\maketitle

\section{Introduction}
In theoretical statistics the marginal polytope plays an important
role. It is the polytope of possible values that a sufficient
statistics can take. It encodes in its face lattice the combinatorial
structure of the boundary of the exponential family defined by the
statistics. For a model on discrete random variables it
can be represented with vertices that have only components 0 or 1,
commonly called a 0/1 polytope.

In coding theory when decoding binary linear codes one can apply
techniques from linear programming and optimize a linear function over
the convex hull of the code words, known as the codeword polytope
\cite{feldman03:_using_linear_progr_to_decod_linear_codes}.

\comment{Our starting point is the fact that for certain choices of
  statistical models on binary random variables the two notions
  coincide. We explore this connection by giving a classification of
  the marginal polytopes that are codeword polytopes.}

Observing that for certain choices of sufficient statistics on binary
random variables these two notions coincide, our main contribution is
a characterization of the corresponding polytopes. We do not address
problems that are directly linked to coding theory. However, we do
hope that our result will contribute to a better understanding of the
closure of exponential families, which is an important problem in
statistics.

The paper is organized as follows: In Section \ref{sec:prelim}, we
introduce the necessary notions to define hierarchical models and fix
the notation. We review different descriptions of so called
interaction spaces in Section
\ref{sec:InteractionSpacesGeneratingSets}. In Section
\ref{sec:CodingTheory}, we establish the link to coding theory.
Finally, in Section \ref{sec:Classification} we give our main result,
the classification of all such full dimensional polytopes whose
vertices form a linear code and give a recursive formula for their
number.

\section{Preliminaries} \label{sec:prelim}
\subsection{Exponential Families of Hierarchical Models}
\label{sec:probabilityprelim} Given a non-empty finite set ${\mathcal
  X}$, we denote the set of probability distributions on ${\mathcal
  X}$ by $\overline{\mathcal P}({\mathcal X})$.  The \emph{support} of $p
\in \overline{\mathcal P}({\mathcal X})$ is defined as $\supp(p) :=
\{x \in {\mathcal X} : p(x) > 0 \}$.
The set of distributions with full support is denoted as
$\mathcal{P}(\mathcal{X})$.
The set $\overline{\mathcal{P}}(\mathcal{X})$ has the geometrical structure of a
$(\abs{\mathcal{X}} - 1)$-dimensional simplex lying in an affine hyperplane of 
\[
\mathbb{R}^{\mathcal{X}} \defas \set{f : \mathcal{X} \to \mathbb{R}},
\]
the vector space of real valued functions on $\mathcal{X}$. Statistical models,
such as hierarchical models are subsets of $\overline{\mathcal{P}}(\mathcal{X})$. In this
paper, we will only consider so called exponential families which are smooth manifolds. 
\begin{defn}
The map
\[
   \exp: \; {\mathbb R}^{\mathcal X} \to {\mathcal P}({\mathcal X}),\qquad
   f \mapsto 
   \frac{\e^{f}}{\sum_{x \in {\mathcal X}} \e^{f(x)}}, 
\]
is called the \emph{exponential map}. It acts component wise by exponentiating
and normalizing. 
Then, an \emph{exponential family} (\emph{in $\mathcal{P}(\mathcal{X})$}) is
defined as the image $\exp({\mathcal I})$ of a linear subspace
${\mathcal I}$ of ${\mathbb R}^{\mathcal X}$.
\end{defn}
An exponential family $\mathcal{E}$ naturally has full support and is
therefore contained in the open simplex
$\mathcal{P}(\mathcal{X})$. However, to get probability distributions
with reduced support one has to pass to the closure
$\overline{\mathcal{E}}$ with respect to the standard topology of
$\mathbb{R}^{\mathcal{X}}$.\comment{As $\mathcal{X}$ is finite, one
  can use the familiar notions of closure in
  $\mathbb{R}^{\mathcal{X}}$.}

Now we consider a compositional structure of ${\mathcal X}$ induced by the set 
$[N]\defas \set{1,\ldots,N}$. Given a subset $A
\subseteq [N]$, we define 
\[
\mathcal{X}_A \defas \set{0,1}^{A},
\]
\comment{We will also write $\mathcal{X} \defas
\mathcal{X}_{[N]}$.} and the natural projection
\[
X_A : {\mathcal X_{[N]}} \to {\mathcal X}_A, \qquad (x_i)_{i \in [N]}
  \mapsto (x_i)_{i \in A} \; .
\]
In the following, we will abbreviate $\mathcal{X} \defas
\mathcal{X}_{[N]}$.  One can view
$\overline{\mathcal{P}}(\mathcal{X})$ as the set of joint probability
distributions of the binary random variables $\set{X_i : i \in [N]}$.
We now use the compositional structure of ${\mathcal X}$ in order to
define exponential families in ${\mathcal P}({\mathcal X})$ given by
interaction spaces. Now, decompose $x\in {\mathcal X}$ in the form $x
= (x_A,x_{[N]\setminus A})$ with $x_A\in {\mathcal X}_A$, $x_{[N]
  \setminus A}\in{\mathcal X}_{[N] \setminus A}$, and define
${\mathcal I}_A$ to be the subspace of functions that do not depend on
the configurations $x_{[N] \setminus A}$:
\begin{eqnarray*}
  {\mathcal I}_A & := & 
  \left\{ f \in {\mathbb R}^{\mathcal X} \; : \; f(x_A,x_{[N] \setminus A}) =  
    f(x_A, x_{[N] \setminus A}')\right. \\
  &    & \qquad \qquad \quad \left. \mbox{for all
      $x_A \in {\mathcal X}_A$, and all $x_{[N] \setminus A}, x_{[N] \setminus A}'
      \in {\mathcal X}_{[N] \setminus A}$}\right\} \; .
\end{eqnarray*}
In the following, we apply these interaction spaces as building blocks
for more general interaction spaces and associated exponential
families \cite{darrochspeed83}. 
The definition of a hierarchical model is based on the notion of a 
hypergraph \cite{lauritzen96}:
\begin{defn}
  A \emph{pre-hypergraph} $\mathcal{A}$ is a non-empty subset of $2^{[N]}
  \setminus{\set{\emptyset}}$ that contains all atoms $\set{i}$ for $i\in [N]$.

  A \emph{hypergraph} is a pre-hypergraph that is (inclusion) complete
  in the following sense: If $A\in \mathcal{A}$ and $\emptyset \neq
  B\subset A$ it follows that $B\in \mathcal{A}$.
\end{defn}

\begin{rem}
  For technical convenience, we have defined hypergraphs to be
  complete.  In this way, it is easy to define a hierarchical model
  for each hypergraph. However, the notion of a pre-hypergraph turns
  out to be more natural in the context of the polytopes and linear
  codes that we consider below.
\end{rem}

Given a hypergraph, we define the associated interaction space by
\[ {\mathcal I}_{\mathcal A} \; := \; \sum_{A \in {\mathcal A}}
{\mathcal I}_A. \] Note that, since a function that depends only on
its arguments in $A$, only depends on its arguments in $B\supset A$,
it suffices to consider the inclusion maximal elements in
$\mathcal{A}$. We denote them by $\mathcal{A}^m$ and have
\[
\mathcal{I}_\mathcal{A} = \sum_{A\in\mathcal{A}^m} \mathcal{I}_A.
\]
We consider the corresponding exponential family:
\begin{defn}
  The \emph{hierarchical model} assigned to the hypergraph
  $\mathcal{A}$ is the exponential family
\[ 
{\mathcal
  E}_{\mathcal A} :=
\exp({\mathcal I}_{\mathcal A}).
\]
\end{defn}
We give two examples for hypergraphs:
\begin{exmp} \label{graphmod}$ $\\
  {\bf (1) Graphical models:} Let $G = (V,E)$ be an undirected graph,
  and define
\begin{equation*}  
  {\mathcal A}_G \defas \{ \emptyset \neq C \subseteq V : 
  \mbox{$C$ is a clique with respect to $G$}\} \, .
\end{equation*} 
Here, a {\em clique\/} is a set $C$ that satisfies the following
property:
\[
a,b \in C, \;\; a \not= b \quad \Rightarrow \quad \mbox{
$\set{a,b} \in E$} \, .
\]
The exponential family ${\mathcal E}_{{\mathcal A}_G}$ is
characterized by Markov
properties with respect to $G$ (see \cite{lauritzen96}). \\
{\bf (2) Interaction order:} The hypergraph associated with a given
interaction order $k \in \{1,2, \dots,N\}$ is defined as
\begin{equation*}
  \akn \defas \set{ \emptyset \neq A \subset [N] : \abs{A} \leq k}.
\end{equation*} 
If appropriate, we will sometimes drop the $N$ and write
$\mathcal{A}_k$. We have defined a corresponding hierarchy of
exponential families studied in \cite{Amari01,ayknauf06}:
\[
    {\mathcal E}_{{\mathcal A}_1} 
    \subseteq {\mathcal E}_{{\mathcal A}_2}  \subseteq 
    \dots \subseteq {\mathcal E}_{{\mathcal A}_N} = 
    {\mathcal P}({\mathcal X}). 
\] 
\end{exmp}
The elements of this hierarchy have nice interpretations. It can be seen that the
closure of the family $\mathcal{E}_{\mathcal{A}_1}$ contains exactly all probability distributions that
factor. This means that 
\begin{equation*}
  p \in \overline{\mathcal{E}_{\mathcal{A}_1}}  \Leftrightarrow p = \prod_{i\in
  [N]} p_i(x_i), 
\end{equation*} where $p_i(x_i)$ are the marginal distributions of $p$.
Generally, an element $p\in
\mathcal{E}_k$ will allow a factorization as 
\begin{equation*}
  p = \prod_{A\subset [N] , \abs{A} = k} \phi_A(x_A),
\end{equation*} 
where $\phi_A$ depends only on $x_A$. However, the $\phi$ are not
necessarily probability distributions and not unique. Note that $p \in
\overline{\mathcal{E}_k}\setminus\mathcal{E}_{k}, k\geq 2$, does not
necessarily admit such a product structure.

We will clarify these definitions in the following simple
\begin{exmp}
Consider the case $N=2$. The configuration space is given as 
\[
\mathcal{X} \defas \set{0,1}^2 = \set{(0,0),(0,1),(1,0),(1,1)}.
\]
The vector space of real valued functions $\mathbb{R}^{\mathcal{X}}$ is
4-dimensional and the probability measures form a 3-dimensional tetrahedron. 
Considering the hypergraphs of fixed interaction order and their exponential
families, one has only two examples here: $\mathcal{E}_{\mathcal{A}_{1,2}}$ and
$\mathcal{E}_{\mathcal{A}_{2,2}} = \mathcal{P}(\mathcal{X})$,
only the first being nontrivial. 
\begin{figure}[htbp]
  \centering
\begin{center}
\setlength{\unitlength}{1cm}
\begin{picture}(7,7)
\put(0,0){\includegraphics{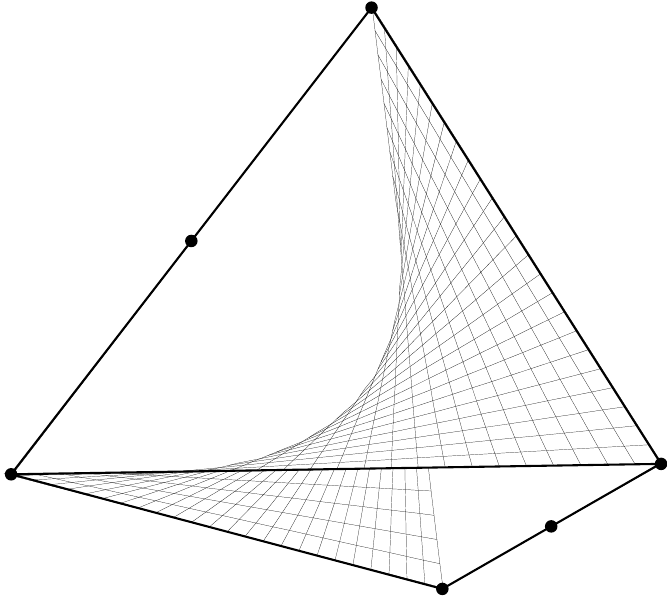}}
\put(3.4,6.28){\footnotesize $\delta_{(1,1)}$}
\put(-0.75,3.7){\footnotesize 
$\frac{1}{2}\left(\delta_{(0,0)} + \delta_{(1,1)}\right)$}
\put(5.5,0.25){\footnotesize 
$\frac{1}{2}\left(\delta_{(1,0)} + \delta_{(0,1)}\right)$}
\put(6.9,1.25){\footnotesize $\delta_{(1,0)}$}
\put(4.1,-0.23){\footnotesize $\delta_{(0, 1)}$}
\put(-0.8,1.2){\footnotesize $\delta_{(0,0)}$}
\put(5.5,1.6){$\boldsymbol{\mathcal E_{\mathcal{A}_{1,2}}}$}
\end{picture}
\end{center}
\caption{The exponential family $\mathcal{E}_{\mathcal{A}_{1,2}}$ in the simplex of 
probability distributions.} \label{fig:simplex}
\end{figure}
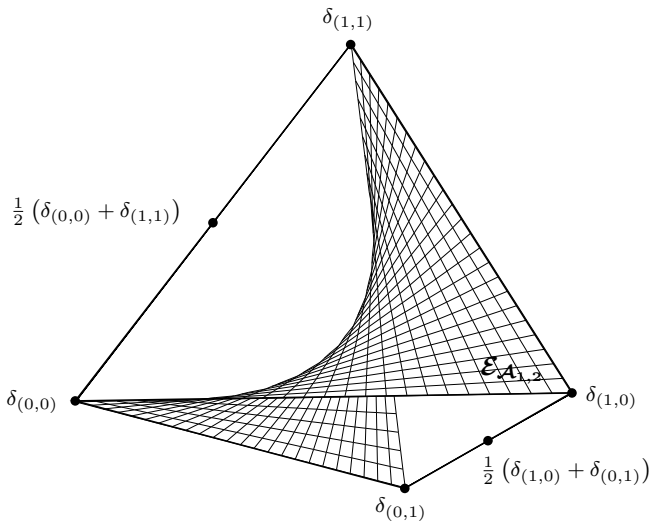
Figure \ref{fig:simplex} shows the situation. The exponential family
$\mathcal{E}_{1,2}$ is a two-dimensional manifold lying inside the simplex. One
should already think about this as a square (the two dimensional cube) molded into
the simplex.
\end{exmp}
In the following, we will study the interaction spaces more thoroughly
by comparing different generating systems. 

\section{Generating Systems of Interaction Spaces}
\label{sec:InteractionSpacesGeneratingSets}
In this section, let $\mathcal{A}$ be fixed. In statistics, different
representations of exponential families have been considered, each of
which has its own benefits and highlights different aspects. We will
review a number of these representations. The particular choice of
parity functions will allow us to make a link to coding
theory. \comment{then add our own which, at least to the authors
  knowledge has not been considered so far.}

As we have introduced exponential families, the key concept is the
interaction space $\mathcal{I}_{\mathcal{A}}$ which is sometimes also called
tangent space to the exponential family. This space completely
characterizes the exponential family. However, there is a choice of
the parameterization of this space, which has been made differently in
different fields. Speaking in terms of linear algebra, one has to
choose a generating system of a linear space.

Let $\mathfrak{B} \defas \set{b_{k} : k\in K}$ be any finite
generating system of $\mathcal{I}_{\mathcal{A}}$. Each such choice gives a
different parameterization of the exponential family and a different
sufficient statistics. The parameterization is identifiable if
$\mathfrak{B}$ is a basis. The exponential family is parameterized as 
\begin{equation*}
  \mathcal{E}_{\mathcal{A},\mathfrak{B}} =
  \set{p\in\mathcal{P}(\mathcal{X}) : p(x) = Z^{-1}_{\theta} \exp
    \left(
      \sum_{k\in K} \theta_{k} b_{k}(x)
    \right)
  : \theta \in \mathbb{R}^{d}},
\end{equation*}
where again $Z_{\theta}$ is the normalization and $d =
\abs{\mathfrak{B}}$ equals the number of parameters. In statistical
physics the exponent is commonly called the energy.

To each choice of $\mathfrak{B}$ there is a polytope constructed as
follows. Consider the vectors
\begin{equation*}
  b_{K}(x) = (b_{k}(x))_{k\in K}
\end{equation*}
Each such vector has as its components the evaluation of every element
in $\mathfrak{B}$ at $x$. The polytope is 
\begin{equation*}
  P_{\mathfrak{B}} \defas \conv \set{ b_{K}(x) : x \in \mathcal{X}}.
\end{equation*}
Since $\mathcal{A}$ contains all atoms, it can be seen that the
polytope has $\abs{\mathcal{X}}$ vertices and the dimension equals the
dimension of the exponential family. By applying some classical
theorems from statistics, such as the existence and uniqueness of
maximum likelihood estimates
\cite{kullback68:_infor_theor_and_statis,csiszar75:_i_diver_geomet_of_probab},
it can be seen that the points of the polytope are in one to one
correspondence with points in the closure of the exponential
family. As we have introduced it here, it is clear that the different
choices of $\mathfrak{B}$ yield different representation of the same
polytope in the sense that they are all affinely equivalent. In
particular, they have the same face lattice.

The polytope $P_{\mathfrak{B}}$ encodes in its face lattice the
combinatorial structure of the exponential family in the sense that a
knowledge of the face lattice gives precise knowledge about the
supports of elements in the closure of the exponential
family. However, direct computation is infeasible for real world
problems.

In statistical physics, and also for various inference methods it is
of interest to compute the free energy, given as the logarithm of the
partition function.  There, variational principles and the techniques
of Legendre transform are applied.  In this setting the points in the
polytope are then the so called dual parameters. See for instance
\cite{wainwright03:_variat_infer_in_graph_model}.

We will review a number of choices for $\mathfrak{B}$:

\paragraph{Statistical Physics - Potentials}
\label{sec:stat-phys-potent}
In statistical physics one considers so called
potentials \cite{winkler_image_analysis,georgii88}. A potential is a
collection of functions $U_{A}, A\subset N$, where $U_{A} \in
\mathcal{I}_{A}$ and $U_{\emptyset} = 0$, such that the energy can be
written as a linear combination hereof. Typically one has a
distinguished state $o$ called the vacuum. A potential is called
\emph{normalized} if $U_{A}(x) = 0$ as soon $x_{i}=o_{i}$ for some
$i\in A$. Given a strictly positive distribution, a corresponding
normalized potential exists and is unique. In our binary setting,
choosing $(0,0,\ldots,0)$ as the vacuum state, the normalized
potential is given by the functions $U_{A} = c_{A} \prod_{i\in A}
x_{i}$, where $c_{A} \in \mathbb{R}$.

One has $\mathfrak{B} = \set{\prod_{i\in A} x_{i} : A \in
  \mathcal{A}}$, and a basis of the interaction space is given by
$\mathfrak{B}$ together with the constant function $x\mapsto 1$.
Expanding a function $H\in \mathbb{R}^{\mathcal{X}}$ in terms of this
basis was called the $\chi$-expansion in the works of Caianiello
\cite{caianiello86:_neuron_equat_revis_and_compl_solved,caianiello75:_synth_of_boolean_nets_and}.

In the case of pair interactions where the hypergraph is given by
$\mathcal{A}_{2,N}$, the polytope $P_{\mathfrak{B}}$ coincides with
the so called \emph{correlation polytope} \cite{dezalaurent97}.
Extending the terminology to an arbitrary hypergraph $\mathcal{A}$, we
call $P_{\mathfrak{B}}$ the \emph{moment polytope}, as each point in
it is the vector of moments of some distribution.

\paragraph{Marginals}
\label{sec:marg-chi-repr}
One representation of an exponential family is given via the linear
map that computes the marginals. Denote $\mathcal{A}^{m}$ the set of
inclusion maximal sets in $\mathcal{A}$. Consider the linear map
\begin{align*}
  \pi_{\mathcal{A}} : \mathbb{R}^{\mathcal{X}} & \to \bigoplus_{A \in
    \mathcal{A}^{m}} \mathbb{R}^{\mathcal{X}_A} \\
  u & \mapsto \left(u_A\right)_{A \in \mathcal{A}}.
\end{align*}
That, for a given vector $u$ computes the set of its maximal marginals
defined as
\begin{equation*}
  u_A(x_A) \defas \sum_
  {y : X_{A}(y) = x_{A} }  u(y).
\end{equation*}
When represented as a matrix with respect to the canonical basis,
$\pi_{\mathcal{A}}$ has rows indexed by pairs $(A,y_{A})$ of a set
$A\in\mathcal{A}^{m}$ and a configuration
$y_{A}\in\mathcal{X}_{A}$. The columns are indexed by configurations
$x \in \mathcal{X}$. Each component then contains the value of the
indicator $\Ind{X_{A}=y_{A}}$.

Denote the $x$-th column of this matrix as $\pi_{x}$ then, the
exponential family is parameterized as
\begin{equation*}
  p(x) = Z_{\theta}^{-1} \exp (
  \left\langle \theta , \pi_{x} \right\rangle) , 
  \qquad \theta\in\mathbb{R}^{d}.
\end{equation*}

In terms of these vectors, the polytope is commonly called the marginal
polytope. It is represented as 0/1 polytope embedded in a high
dimensional space.

\paragraph{An orthogonal basis of characters}
\label{sec:basis-characters_desc}
In the binary case $\mathcal{X} = \set{0,1}^{[N]}$, a natural basis
for $\mathbb{R}^{\mathcal{X}}$ is given by the characters of
$\mathcal{X}$. Here, we assume pointwise addition modulo 2 as the
group operation. For every subset $A\in\mathcal{A}$ define the
function $e_A : \mathcal{X} \to \set{-1,1}$ by
\begin{align*}
  e_A(x) \defas (-1)^{E(A,x)}
\end{align*}
where $E(A,x) \defas \abs{\set{i\in A : x_i = 1}}$.  It can be seen
that, if $\mathcal{A}$ is a hypergraph, $\set{e_A : A \in
  \mathcal{A}}$ together with the constant function $e_{\emptyset} :x\mapsto 1$ is an orthogonal basis of the interaction space
$\mathcal{I}_{\mathcal{A}}$. This approach was followed in
\cite{kahleay06}. Various people, starting with
Caianiello \cite{caianiello75:_synth_of_boolean_nets_and} have called
this the $\eta$-expansion. Note that if one considers random variables
taking values in $\set{\pm 1}$ this basis equals the monomial basis
$\set{\prod_{i\in A} x_{i} : A\subset [N]}$ considered above.

\paragraph{A basis of parity functions}
\label{sec:basis-parity-funct}

Finally, we will introduce yet another basis of
$\mathcal{I}_{\mathcal{A}}$ which is derived from the basis of
characters. To each $\emptyset \neq A \subset [N]$, we define a vector
in $\mathbb{R}^{\mathcal{X}}$.
\begin{equation}
  \label{eq:fAdefinition}
  f_A (x) \defas 
  \begin{cases}
    1 & \text{ if } \abs{\supp(x) \cap A} \text{ is odd} \\
    0 & \text{ otherwise.}
  \end{cases}
\end{equation}
The following proposition is easily checked:
\begin{prop}
  \label{sec:fa_basis_prop} \sloppy Let $\mathbbm{1}: \mathcal{X}
  \to \mathbb{R}$ be the constant function $x \mapsto 1$.  The set
  $\set{f_A : A \in \mathcal{A}} \cup \set{\mathbbm{1}}$ is a basis of
  $\mathcal{I}_\mathcal{A}$.
\end{prop}

One crucial point about choosing this representation is that it gives,
if the constant function is omitted, a full dimensional 0/1 polytope,
the vertices of which form an additive group and thereby a linear code
(see Proposition \ref{sec:parity-code}). For all other choices of
$\mathfrak{B}$ discussed in this section the image of $b_{K}$ is not a
subgroup of $\set{0,1}^{d}$ or the multiplicative group $\set{\pm
  1}^{d}$.

While in the construction of a hierarchical model we assumed a
hypergraph, the following polytope is an interesting object of study
also in the general case of a pre-hypergraph:
\begin{defn}
  Let $\mathcal{A}$ be a pre-hypergraph. We define
  \begin{equation*}
    \mathcal{F}_{\mathcal{A}} \defas \conv \set{ f_{\mathcal A} (x) :
      x \in \mathcal{X}}.
  \end{equation*}
\end{defn}

If $\mathcal{A}$ is a hypergraph, then this is affinely equivalent to
the marginal polytope of the corresponding exponential family.  In the
case of the hypergraphs $\akn$ we write $\fkn \defas
\mathcal{F}_{\akn}$. The rest of the paper is devoted to the study of
this class of polytopes.

\begin{rem}[CUT-Polytopes] 
  There is a well known \cite{dezalaurent97} affine equivalence
  between CUT polytopes of graphs \cite{ziegler00LecturesPolytopes} and
  binary marginal polytopes: 

  Namely, to each graph $G$ we can associate the hypergraph
  $\mathcal{A}_{G} = V(G) \cup E(G)$. This is distinct from what was
  called a graphical model above, as not the cliques are
  considered. Some authors refer to the corresponding statistical
  model as a \emph{graph model}. From $G$ construct the coned graph
  $\hat{G}$ with an additional vertex:
  \begin{equation*}
    V(\hat{G}) \defas V(G) \cup \set{*},
  \end{equation*}
  and edges 
  \begin{equation*}
    E(\hat{G}) \defas E(G) \cup \set{(v,*) : v \in V(G)}.
  \end{equation*}  Then, 
  denoting the CUT polytope of $\hat{G}$ as $CUT(\hat{G})$ one has
  \begin{equation*}
    \mathcal{F}_{\mathcal{A}_{G}} = CUT(\hat{G}).
  \end{equation*}
  Using the representation in terms of the vectors
  $f_{\mathcal{A}_{G}}(x)$, $x\in \mathcal{X}$, the proof of this
  equivalence becomes a simple renaming of coordinates.
\end{rem}

\begin{rem}[Covariance Mapping]
  As remarked above, in the representation with monomials $\prod_{i\in
    A} x_{i}$ one finds the correlation polytope $COR(N)$ as a special
  case. From the last remark it follows that the CUT-polytope of the
  complete graph $K_{N+1}$ is equal to
  $\mathcal{F}_{\mathcal{A}_{2,N}}$. There exists an affine
  equivalence between $COR(N)$ and $CUT(K_{N+1})$ called the
  \emph{covariance mapping} \cite{dezalaurent97}. It can be seen that
  this mapping generalizes to a mapping between binary marginal
  polytopes and the corresponding moment polytopes. It therefore
  might be suitable to consider the parity representations
  $\mathcal{F}_{\mathcal{A}}$ of binary marginal polytopes for a
  generalization of CUT-polytopes to arbitrary (pre)-hypergraphs.
\end{rem}

\subsection{Computations and elementary properties}
\label{sec:comp-elem-prop}

Using the geometry software \texttt{polymake}
\cite{GawrilowJoswigPolymake}, one can compute linear descriptions of
polytopes. As an example, we give here the F-Vectors of
$\mathcal{F}_{k,N}$ for the cases $N=3,4$.  For $N=5$, the F-Vector is
too complicated to be computed by the brute force approach of
\texttt{polymake}. However, waiting sufficiently long, one can get the 6800
facet defining inequalities of $\mathcal{F}_{3,5}$ and the 3835488
facets of $\mathcal{F}_{4,6}$.

\begin{exmp}
  In Tables \ref{tab:3node} and \ref{tab:4node}, we give the F-Vectors
  of $\mathcal{F}_{k,N}$ for $N=3,4$, computed using
  \texttt{polymake}. The rows label the dimension of the faces, the
  columns the value of $k$. The reader might wonder about the fact
  that the face lattices of $\fkn$ are up to a certain dimension
  isomorphic the face lattice of the simplex. This property, commonly
  called neighborliness, follows from a general result in
  \cite{kahle08_degree}. The last row refers to whether the polytope
  is simple or not.
\begin{table}[htpb]
\centering
\begin{tabular}{|l|c|c|c|}
  \hline
  $d \backslash k $ & 1 & 2 & 3 \\
  \hline
  0 & 8 & 8 & 8 \\
  1 & 12 & 28 & 28 \\
  2 & 6 & 56 & 56 \\
  3 & 1 & 68 & 70 \\
  4 & - & 48 & 56 \\
  5 & - & 16 & 28 \\
  6 & - & 1 & 8 \\
  7 & - & - & 1 \\
  sum & 27& 225 & 255\\
  simple & y& n & y \\ 
  \hline
\end{tabular} 
\caption{Face structure of $\mathcal{F}_{k,3}$}
\label{tab:3node}
\end{table}

\begin{table}[htpb]
  \centering
  \begin{tabular}{|l|c|c|c|c|}
    \hline
    $d \backslash k$ & 1 & 2 & 3 & 4 \\
    \hline
    0 & 16 & 16 & 16 & 16 \\
    1 & 32 & 120 & 120 & 120 \\
    2 & 24 & 560 & 560 & 560 \\
    3 & 8 & 1780 & 1820 & 1820 \\
    4 & 1 & 3872 & 4368 & 4368 \\
    5 & - & 5592 & 8008 & 8008 \\
    6 & - & 5060 & 11440 & 11440  \\
    7 & - & 2600 & 12868 & 12870 \\
    8 & - & 640 & 11424 & 11440 \\
    9 & - & 56 & 7952 & 8008 \\
    10 & - & 1 & 4256 & 4368 \\
    11 & - & - & 1680 & 1820 \\
    12 & - & - & 448 & 560 \\
    13 & - & - & 64 & 120 \\
    14 & - & - & 1 & 16 \\
    15 & - & - & - & 1 \\
    sum & 81 & 20297 & 65025 & 65535\\
    simple & y& n & n & y \\ 
    \hline
  \end{tabular}  
  \caption{Face structure of $\mathcal{F}_{k,4}$}
  \label{tab:4node}
\end{table}
\end{exmp}

In the following, we will list elementary properties of $\fkn$ that follow easily
from the definition. 
\begin{enumerate}[(i)]
\item $\mathcal{F}_{1,N}$ is the $N$-cube.
\item $\mathcal{F}_{N,N}$ is the $(2^{N}-1)$-dimensional simplex.
\item every $\mathcal{F}_{k,N}$ has dimension $d =
  \abs{\akn}$.
\item every $\fkn$ has $2^{N}$ vertices.
\item $(0,\ldots,0)$ is a vertex.
\item every $\fkn$ is a projection of the $(2^N-1)$-dimensional
  simplex $\mathcal{F}_{N,N}$ along coordinate axes.
\item For every $\fkn$, there is a projection along coordinate
  axes that projects it to the $N$-cube $\mathcal{F}_{1,N}$.
\end{enumerate}

\begin{rem}
  In \cite{hostensullivant02a} it was remarked that
  $\mathcal{F}_{N-1,N}$ has exactly $4^{N-1}$ facets. The extreme
  points of these facets are also known.  A set $\mathcal{Y}\subsetneq
  \mathcal{X}$ defines a face if and only if it contains neither
  $\mathcal{U} \defas \set{x \in \mathcal{X} : f_{[N]}(x) = 1}$ nor
  its complement. Note that the set $\mathcal{U}$ and its complement
  are exactly the set of configurations with a fixed parity. As the
  vertices of $\mathcal{F}_{N-1,N}$ have only one affine dependency,
  it is not difficult to prove this fact using the Gale transform. By
  the above $\mathcal{F}_{N-1,N}$ is combinatorially isomorphic to the
  so called cyclic polytope \cite{ziegler94}.
\end{rem}

In the following, we develop the connection to coding theory.
\section{A Link to Coding Theory}
\label{sec:CodingTheory}
We briefly recall the definition of a
linear code. For a detailed introduction into coding theory see
for instance \cite{vanLint99}.
\label{sec:codingprelim}
Consider the finite field $\mathbb{F}_2 = (\set{0,1}, \oplus, \odot)$ with
addition and multiplication mod 2. In coding theory, one studies particularly vector
spaces over this field.
\begin{defn}
  A binary \emph{$[n,k]$-linear code} is a linear subspace $L$ of
  $\mathbb{F}_2^n$ such that $\dim L = k$. A \emph{generator matrix}
  $G$ for $L$ is a $k$ by $n$ matrix which has as its rows a basis of
  $L$. Given $L$ one can find an equivalent\footnote{Two codes are called
  equivalent if one can be transformed into the other by applying a permutation
  on the positions in the codewords, and for each position a permutation of the
  symbols.} code such that it has a generator matrix in
  \emph{standard form,} \ie $G=(E_k,H)$, where $E_k$ is the
  $k$ by $k$ identity matrix.
\end{defn}

The following proposition states that the vertices of
$\mathcal{F}_{\mathcal{A}}$ form a linear code for any pre-hypergraph
$\mathcal{A}$. A special case of this connection has been mentioned in
Example 2 in \cite{wainwright03:_variat_infer_in_graph_model}.

\begin{prop}
\label{sec:parity-code}
Let $\set{0,1}^{\mathcal{A}}$ be considered as a vector space over
the finite field $\mathbb{F}_2$. Then the image of $\mathcal{X}$ under
$f_{\mathcal{A}}$ is a linear subspace. If we also consider
$\mathcal{X} = \mathbb{F}_2^N$ as a vector space over $\mathbb{F}_2$,
$f_\mathcal{A}$ is an injective homomorphism between vector
spaces. Its image forms an $[\abs{\mathcal{A}},N]$-linear code. A
generator matrix in standard form has as its rows the vectors
$f_{\mathcal{A}}(e_i)$ for $i=1,\ldots,N$, where $e_i$ is the $i$-th
unit vector in $\mathbb{F}^{N}_{2}$.
\end{prop}
\begin{proof}
  Since scalar multiplication is trivial,
  we only need to show 
  \begin{equation}
    \label{eq:linearityofeA}
    f_\mathcal{A}(x\oplus y) = f_\mathcal{A}(x) \oplus f_\mathcal{A}(y) \qquad
    \text{ for $x,y\in\mathcal{X}$}.
  \end{equation}
  Let $A\in\mathcal{A}$, it suffices to show the identity for $f_A$. To do so,
  introduce 
  \begin{align*}
    M & \defas \set{i\in A: (x_i = 1 \land y_i = 0) \lor (x_i=0 \land y_i = 1)}, \\
    M_x & \defas \set{i\in A: x_i = 1}, \\
    M_y & \defas \set{i\in A: y_i = 1}. \\
  \end{align*}
  Then $f_A(x\oplus y) = \abs{M}$, $f_A(x) = \abs{M_x}$, and $f_A(y) =
  \abs{M_y}$.
  We find that $M$ is the symmetric difference of $M_x$ and $M_y$:
\begin{equation*}
  M = M_x \bigtriangleup M_y.
\end{equation*}
Since $\abs{M_x \bigtriangleup M_y} = \abs{M_x} + \abs{M_y} - 2
\abs{M_x \cap M_y}$, we have that in $\mathbb{F}_2$
\begin{equation*}
  \abs{M} = \abs{M_x} \oplus \abs{M_y},
\end{equation*} 
and therefore \eqref{eq:linearityofeA} holds. We now show that
$f_\mathcal{A}$ is injective. To see this, assume that $f_{\mathcal{A}}(x) =
f_{\mathcal{A}}(y)$. Since $\mathcal{A}$ contains all atoms
$\set{i}\subset [N]$, we get for every $i\in [N]$: $f_{\set{i}}(x) =
f_{\set{i}}(y)$. This implies $x_i = y_i$ and, hence, $x=y$. 
Since $\mathcal{X}$ considered as an $\mathbb{F}_2$ vector space has
dimension $N$, also $f_\mathcal{A}(\mathcal{X})$ has dimension $N$ and
therefore forms an $[\mathcal{A},N]$-linear code. 
\end{proof}
\begin{rem}
To write down the generator matrix, one has to impose a numbering on the elements
in $\mathcal{A}$. If the numbering is in such a way that $A_i = \set{i}$ for
$1\leq i \leq N$, then the generator matrix is in standard form.
\end{rem}

\begin{rem} An important property of a linear code is its distance, which is
  defined as the minimal Hamming distance between different elements of the code. 
  For the hierarchical model of the hypergraph $\akn$, the distance of
  the code is given by
  \begin{equation*}
    d = \sum_{l=0}^{k-1}\binom{N-1}{l}.
  \end{equation*}
\end{rem}
\begin{proof}
  Let $d(x,y)$ denote the hamming distance of $x,y\in\mathcal{X}$. If $d(x,y)=1$, 
  then $d(f_{\mathcal{A}_k}(x),f_{\mathcal{A}_k}(y))$ equals the
  number of subsets of $[N]$ which contain a given element and have cardinality at
  most $k$. 
\end{proof}

In the following, we will elaborate on the opposite direction. Let 
$2^N \geq s \geq N$. Assume we are given an $[s,N]$ linear code. Without loss of
generality, we assume that it has a generator matrix in standard form. We will
construct a pre-hypergraph $\mathcal{A}$ from the columns of the generator
matrix. Since $\mathcal{A}$ is a set, while the columns are a list, repetitions
of columns will be lost. If one considers only non-repetitive codes, then our
construction is injective, and the codewords are given by the vertices of
$\mathcal{F}_\mathcal{A}$.

Let $E_N \in \mathbb{F}_2^{N\times N}$ denote the
identity matrix in dimension $N$. Assume the generator matrix $G = (E_N,H)\in
  \mathbb{F}_2^{N\times s}$ has no 2 identical columns. (This implies
  $s\leq 2^{N}$.) Denote by $\set{e_i : i=1,\ldots,N}$ the
  canonical basis of $\mathbb{F}_2^N$. Using the columns of $H$, we define sets
  \begin{equation*}
    A_j \defas \set{i \in {[N]} : H_{ij} = 1}, \qquad j=1,\ldots, s-N
  \end{equation*}
  and then, 
  \begin{equation*}
    \mathcal{A} \defas \set{\set{1},\ldots \set{N}, A_1, \ldots, A_{s-N}}.
  \end{equation*} 
  Note that the elements of $\mathcal{A}$ are numbered in a natural way
  such that we can use $\mathcal{A}$ as an index set for the columns of
  $G=(G_{i,A})_{i=1,\ldots,N,\;A=\set{1},\ldots,\set{N},A_{1},\ldots,A_{s-N}}$.
  
  To see that $\set{f_\mathcal{A}(e_i) : i = 1,\ldots,N}$ is the set of rows
  of the generator matrix, we evaluate
  \begin{equation*}
    f_A(e_i) = \delta_{\{i\in A\}} = G_{i,A}
  \end{equation*}
  which holds by definition of the $A_j$.

Summarizing, every binary linear code (in standard form) 
corresponds to a pre-hypergraph. However, two codes that differ only in
repetitions of columns in the generator matrix will be mapped to the same
pre-hypergraph. Then, if it is a hypergraph, the linear code is the 
linear code of a hierarchical model.

\section{Classification}
\label{sec:Classification}
As we have seen, the polytopes $\mathcal{F}_{\mathcal{A}}$ are full
dimensional polytopes such that the vertices form a linear code. In
this last section, we classify all polytopes with this property. Then
we investigate which of them can be realized as polytopes of
hierarchical models.  For a convex polytope $P$, let $V(P)$ denote the
vertex set of $P$. For $n \in \mathbb{N}$, put
\[
C_n \defas [0,1]^{n},
\]
\[
W_n \defas \set{0,1}^{n} = V(C_n).
\]
Hence $(W_n,\oplus)$ is an Abelian group that is canonically
isomorphic to $(\mathbb{F}_{2}^{n},\oplus)$. We consider $W_n$ as a
subset of $\mathbb{R}^n$ and write ``$\oplus$'' whenever we mean
addition modulo 2, while ``$+$'' means ordinary addition in
$\mathbb{R}^n$.

In the following, we develop an algorithm that determines - by
induction for every $n\in\mathbb{N}$ - all polytopes $P\subset
\mathbb{R}^n$ with $V(P)\subset W_n$ satisfying the following
conditions:
\begin{enumerate}[(I)]
  \item \label{cond1} $(V(P),\oplus)$ is a subgroup of $(W_n,\oplus)$. 
  \item \label{cond2} $P$ has dimension $n$.
\end{enumerate}

Note that the number of vertices of such a polytope is a power of two. 
Of course, the full $n$-dimensional cube $P = C_n$ satisfies (\ref{cond1}) and
(\ref{cond2}). To start the induction, we remark that there are no further such
polytopes in the cases $n=1$ and $n=2$. For $n=3$, the 3-dimensional regular
simplex $S$ with 
\[
V(S) = \set{(0,0,0),(1,1,0),(1,0,1),(0,1,1)}
\]
satisfies (\ref{cond1}) and (\ref{cond2}), too.

More generally, by \cite[Theorem 2.2]{Wenzel2006}, we have the following
\begin{prop}
  For $n\geq 3$, the following statements are equivalent:
  \begin{enumerate}
    \item $(W_n,\oplus)$ contains some subgroup $U$ such that $\conv (U)$ is a
      regular simplex of dimension $n$.
    \item $n+1$ is some power of 2.
  \end{enumerate}
\end{prop}

In the case $n=3$, the full 3-cube as well as the regular simplex mentioned
above are the only polytopes satisfying conditions (\ref{cond1}) and
(\ref{cond2}). Note that also 
\[
\set{(0,0,0),(1,0,0),(0,1,1),(1,1,1)}
\]
determines a subgroup of $(W_3,\oplus)$; however, the convex closure has
dimension 2.

For fixed $n\geq 2$, define the bijections 
$
\pi_0 : \mathbb{R}^n \times \set{0} \to \mathbb{R}^n$ and $\pi_1:
\mathbb{R}^n\times\set{1} \to \mathbb{R}^n$ by 
\begin{align*}
  \pi_0 (x_1,\ldots,x_n,0) & \defas (x_1,\ldots,x_n), \\
  \pi_1 (x_1,\ldots,x_n,1) & \defas (x_1,\ldots,x_n). \\
\end{align*}
For $1\leq i \leq n$ put 
\begin{align*}
  H_i & \defas \set{(x_1,\ldots,x_n) \in \mathbb{R}^{n} : x_i=0}, \\
  H_i' & \defas \set{(x_1,\ldots,x_n) \in \mathbb{R}^{n} :
  x_i=\frac{1}{2}}.
\end{align*}
Moreover, let $z\defas \left( \frac{1}{2},\ldots, \frac{1}{2} \right)$ denote
the center of the $n$-cube $C_n$.

To determine recursively all 0/1-polytopes $P\subset\mathbb{R}^{n}$ that fulfill
(\ref{cond1}) and (\ref{cond2}), we prove first the following 
\begin{prop}
  \label{sec:main_equiv_thm}
  Suppose that $n\geq 2$ and that $P\subset\mathbb{R}^{n}$ is a 0/1-polytope
  satisfying (\ref{cond1}) and (\ref{cond2}). Assume that $(U,\oplus)$ is a
  subgroup of $(V,\oplus) \defas (V(P),\oplus)$ with $\abs{V:U} =2$. Then the
  following statements are equivalent:
  \begin{enumerate}[(i)]
    \item\label{st1} The polytope $Q \subset \mathbb{R}^{n+1}$, given by 
      \begin{equation}
        \label{eq:liftpolytope}
      V(Q) = \pi_0^{-1}(U)\cup\pi_1^{-1}(V\setminus U),
    \end{equation}
      has dimension $n+1$.
    \item\label{st2} There \textbf{does not} exist some index $i$ with $1\leq i\leq n$ such
      that $U\subset H_i$. In other words, none of the affine hyperplanes $H_i'$
      separates $\conv(U)$ and $\conv(V\setminus U)$.
    \item\label{st3} One has $z \in \conv(U) \cap \conv(V\setminus U)$.
    \item\label{st4} One has $\conv(U) \cap \conv(V\setminus U) \neq \emptyset$.
  \end{enumerate}
\end{prop}
\begin{proof}
  \textbf{(\ref{st1}) $\to$ (\ref{st2}): }\\
  Suppose that $U\subset H_i$ holds for some $i$ with $1\leq i \leq n$. Put 
  \[
  \tilde{H}\defas \set{(x_1,\ldots,x_n,x_{n+1})\in\mathbb{R}^{n+1} :
  x_{n+1} = x_i}.
  \]
  Since $P=\conv(V)$ has dimension $n$ and since $\abs{V:U} = 2$, we must have
  $x_i = 1$ whenever $(x_1,\ldots,x_n) \in V\setminus U$. This means that $V(Q)$
  - and hence also $Q$ - is contained in the $n$-dimensional hyperplane
  $\tilde{H}$, in contradiction to $(\ref{st1})$.

  \textbf{(\ref{st2}) $\to$ (\ref{st3}): }\\
  For $1\leq i \leq n$, let $\alpha_i : \mathbb{F}_2^n \to \mathbb{F}_2 $ denote
  the linear map given by $\alpha_i(x_1\ldots,x_n) \defas x_i$. By assumption,
  $\alpha_{i \upharpoonright U} $ is surjective for $1 \leq i \leq n$. Hence we have 
  \[
  \abs{\set{u\in U : \alpha_i(u) = 0}} = \abs{\set{u \in U : \alpha_i(u) = 1}}
  \quad \text{ for } 1 \leq i \leq n.
  \]
  This means that 
  \[
  z = \frac{1}{\abs{U}} \sum_{u\in U} u \in \conv(U),
  \]
  where the sum is taken in $\mathbb{R}^n$.
  
  Now fix $v_1\in V\setminus U$. Since $V\setminus U = \set{v_1 \oplus u : u\in
  U}$, we get also
  \[
  \abs{\set{v\in V\setminus U  : \alpha_i(v) = 0}} = \abs{\set{v\in V\setminus U
  : \alpha_i(v) = 1}} \quad\text{ for } 1 \leq i \leq n,
  \]
  and hence 
  \[
  z = \frac{1}{\abs{V\setminus U}} \sum_{v \in V\setminus U} v \in \conv
  (V \setminus U).
  \]

  \textbf{(\ref{st3}) $\to$ (\ref{st4}) } is trivial.

  \textbf{(\ref{st4}) $\to$ (\ref{st1}):}\\
  Consider the projection $\pi : \mathbb{R}^{n+1} \to \mathbb{R}^n $ given by
  \[\pi(x_1,\ldots,x_n,x_{n+1}) = (x_1,\ldots,x_n).\] Suppose that the assertion
  is wrong; hence $Q$ is contained in some $n$-dimensional -homogeneous-
  hyperplane $G \subset \mathbb{R}^{n+1}$.
  Since 
  \[
  \pi(V(Q)) = U \dot{\cup} (V\setminus U) = V,
  \]
  the polytope $Q$ has the same dimension as $P=\conv(V)$, that is $n$.
  Thus, the restriction $\pi_{\restriction G}$ is a linear isomorphism from $G$
  onto $\mathbb{R}^n$, and there exists some $\mathbb{R}$-linear map $\alpha:
  \mathbb{R}^n \to G$ satisfying 
  \[
  (\alpha \circ \pi)(w) = w \qquad \text{for all } w \in G.
  \]
  By definition of $V(Q)$, this means:
  \begin{align*}
    \alpha(U) & = \pi_0^{-1}(U), \\
    \alpha(V\setminus U) & = \pi_1^{-1} (V\setminus U).
  \end{align*}
  Hence, $\alpha(\conv(U)) = \conv(\alpha(U))$ and $\alpha(\conv(V \setminus U)) = \conv(\alpha(V\setminus U))$
  are linearly separated by the affine hyperplane 
  \[
  K\defas \set{(x_1,\ldots,x_n,x_{n+1}) \in \mathbb{R}^{n+1} : x_{n+1} =
  \frac{1}{2}}.
  \]
  By (\ref{st4}) this is impossible.
\end{proof}

\begin{exmp}
\label{sec:example_pre_hyp}
  We investigate the statement of Proposition \ref{sec:main_equiv_thm}
  for polytopes corresponding to a pre-hypergraph $\mathcal{A}$. We
  start by considering the matrix which has as its rows the vectors
  $f_{\mathcal{A}}(x)$, where $\mathcal{A} = 2^{[N]} \setminus
  \set{\emptyset}$. The rows are labeled by the binary strings of
  length $N$, that is by $\mathcal{X}$, while the columns are indexed
  by the non-empty subsets of $[N]$. Therefore the rows of this matrix
  are the coordinates of the vertices of the simplex
  $\mathcal{F}_{N,N}$:
\begin{center}
\begin{tabular}{|c||c|c|c|c|c|c|c|c|}
  \hline
  $x$ & $\set{1}$ & $\cdots$ & $\set{N} $ & $\set{1,2}$ & $\cdots$ &
  $\set{1,2,3}$ & $\cdots$ & $[N]$ \\
  \hline\hline
 $(00\ldots0)$ &  $0$ & $\cdots$ & $0$ & $0$ & $\cdots$ & $0$ & $\cdots$ & $0$ \\
 \hline
 $(00\ldots1)$ &  $0$ & $\cdots$ & $1$ & $f_{\set{1,2}}(x)$ & $\cdots$ &
 $f_{\set{1,2,3}}(x)$ & $\cdots $ & $f_{[N]}(x)$ \\
 \hline
 \vdots & \vdots &  $\cdots$ & \vdots & \vdots & $\cdots$ & \vdots
 & $\cdots$ & \vdots \\
 \hline
 $(11\ldots1)$ & $1$ & $\cdots$ & $1$ & $f_{\set{1,2}}(x)$ & $\cdots$ &
 $f_{\set{1,2,3}}(x)$& $\cdots $ & $f_{[N]}(x)$ \\
 \hline
\end{tabular}
\end{center}
We note the following facts:
\begin{itemize}
  \item The columns of this matrix are exactly the $2^{N}-1$ non-zero binary
    strings of length $N$.
  \item There are $2^N-1$ subgroups $U$ of index 2 of the $N$-cube, which correspond to
the columns of the matrix. To define them let a column $A$ be fixed, then put $U \defas
\set{x\in \mathcal{X} : f_A(x) = 0}$. The maps $f_A : \mathcal{X} \to
\set{0,1}$ are exactly the $2^N-1$ surjective homomorphisms having the nontrivial
subgroups as their kernels. 
  \item The vertices of every polytope $\mathcal{F}_{\mathcal{A}}$ are given by
    deleting columns  from this matrix that correspond to sets not in $\mathcal{A}$.
\item In particular, by restriction to the first $N$ columns,
we get the vertices of the $N$-cube $\mathcal{F}_{1,N}$.
\end{itemize}
Now, assume that $P$ is the $N$-cube. We choose a column of the
matrix, corresponding to a subgroup of index 2. There are two
possibilities.  If we choose a column corresponding to an atom, then
$(ii)$ is wrong, the dimension does not grow when adding this column
to the coordinates (as we have doubled a coordinate). If, on the other
hand, we choose a column corresponding to a set $A$ with cardinality
two or more, then we are in the situation of Proposition
\ref{sec:main_equiv_thm}, since $(ii)$ holds. The lift
\eqref{eq:liftpolytope} will be full dimensional, and its vertices are
given by the submatrix with columns $\set{1},\ldots,\set{N},A$.
Continuing from here, choosing another subgroup, the dimension will
grow if and only if it does not correspond to one of the sets
$\set{1},\ldots,\set{N},A$.  Iteratively, the choices narrow down and,
finally, when all columns have been chosen, the polytope $Q$ is a
simplex.
\end{exmp}

We will now formalize this procedure. For a fixed polytope $P$ as in
Proposition \ref{sec:main_equiv_thm}, put
\[
U_i \defas V(P) \cap H_i \qquad \text{for } 1 \leq i \leq n.
\]
Clearly, conditions (\ref{cond1}) and (\ref{cond2}) imply that $\abs{V(P) : U_i}
= 2$ holds whenever $1 \leq i \leq n$.

Based on the equivalence of (\ref{st1}) and (\ref{st2}) in Proposition
\ref{sec:main_equiv_thm}, we are now able to prove that the following algorithm
yields recursively all 0/1-polytopes satisfying (\ref{cond1}) and (\ref{cond2}).

\begin{alg} \label{sec:iterative_algorithm} \ \\
  \textbf{Initialization for $n=1$:}
  \begin{itemize}
  \item $\mathfrak{P}_{1} \defas \set{ [0,1] }$.  
  \end{itemize}
  \textbf{Step $n \to n+1$:} Based on $\mathfrak{P}_{n}$ construct a
  new set $\mathfrak{P}_{n+1}$ consisting of all polytopes $Q$ such
  that there exists $P \in \mathfrak{P}_{n}$ with
  \begin{itemize}
  \item $Q = P\times[0,1]$ or
  \item $Q\subset \mathbb{R}^{n+1}$ with
    \begin{equation}
      \label{eq:projection_condition}
      V(Q) = \pi_0^{-1}(U) \cup \pi_{1}^{-1}(V(P)\setminus U),
    \end{equation}
    where $U$ runs through all subgroups of $(V(P),\oplus)$ with
    $\abs{V(P):U} =2$ and $U\neq U_i$ for $1 \leq i \leq n$.
  \end{itemize}
\end{alg}
\begin{rem}
  Note that in the case $Q = P\times [0,1]$, the number of vertices is
  doubled, while in the other cases the number of vertices of $Q$
  equals the number of vertices of $P$. Furthermore, it is interesting
  to see that the two possible operations commute in the following
  sense. Starting from some cube $W_n$, lifting it to $W_{n+1}$ and
  then choosing a subgroup $U$ to apply the lift
  \eqref{eq:projection_condition} gives the same polytope as choosing
  the subgroup $\pi(U)$ from $W_n$ and then taking the prism over the
  lifted polytope, where $\pi: \mathbb{R}^{n+1} \to \mathbb{R}^{n}$ is
  the canonical projection.  Therefore, all polytopes that are
  constructed by the algorithm can be thought of as lifted cubes
  $W_n$.
\end{rem}

The classification will be complete with:
\begin{thm}
  \label{sec:algothm} 
  For all $n\in\mathbb{N}$, the set $\mathfrak{P}_{n}$ in Algorithm
  \ref{sec:iterative_algorithm} consists of all $n$-dimensional 0/1
  polytopes that satisfy conditions (\ref{cond1}) and (\ref{cond2}).
\end{thm}
\begin{proof}
  First we show that all polytopes $Q \in \mathfrak{P}_{n+1}$ satisfy
  conditions (\ref{cond1}) and (\ref{cond2}), with $n$ replaced by
  $n+1$. This is clear in the case of the prism $Q = P\times [0,1]$.

  If $Q$ satisfies \eqref{eq:projection_condition}, then clearly
  $(V(Q),\oplus)$ is a subgroup of $(W_{n+1},\oplus)$, because $U$ is
  a subgroup of $(V(P),\oplus)$ with $\abs{V(P):U}=2$. Moreover,
  (\ref{st2}) $\to$ (\ref{st1}) in Proposition
  \ref{sec:main_equiv_thm} implies that $Q$ has dimension $n+1$,
  because $U\neq U_i$ for $1 \leq i \leq n$. Hence, $Q$ satisfies
  conditions (\ref{cond1}) and (\ref{cond2}).
  
  Vice versa, assume that $Q\subset \mathbb{R}^{n+1}$ fulfills (\ref{cond1}) and
  (\ref{cond2}). Consider again the projection $\pi : \mathbb{R}^{n+1} \to
  \mathbb{R}^n$ onto the first $n$ coordinates, and put $P\defas \pi(Q)$. Since
  $Q$ has dimension $n+1$, $P$ has dimension $n$. If $\pi_{\restriction
  V(Q)}$ is not injective, then $Q$ is the prism $P \times [0,1]$, because
  $\left( V(Q),\oplus \right)$ is a subgroup of $(W_{n+1},\oplus)$.
  If $\pi_{\restriction V(Q)}$ is injective, put 
  \[
  U \defas \set{(x_1,\ldots,x_n) \in V(P) | (x_1,\ldots,x_n,0) \in Q}.
  \]
  Then $(U,\oplus)$ is a subgroup of $(V(P),\oplus)$ with
  $\abs{V(P):U}=2$, because $Q$ has dimension $n+1$. Moreover,
  equation \eqref{eq:projection_condition} holds for $U$ as just
  defined.  Finally, Proposition \ref{sec:main_equiv_thm}, (\ref{st1})
  $\to$ (\ref{st2}) shows that $U\neq U_i$ for $1 \leq i \leq
  n$. Hence, our algorithm includes the determination of $Q$.
\end{proof}

As a first application of Theorem \ref{sec:algothm} we can count the
number of $n$-dimensional polytopes that satisfy conditions
(\ref{cond1}) and (\ref{cond2}). Let $c_{n} \defas
\abs{\mathfrak{P}_{n}}$. For $1\leq k \leq n$, let $c_{n}(k)$ denote
the number of all 0/1 polytopes $P \subset \mathbb{R}^{n}$ with
$\abs{V(P)} = 2^{k}$ that satisfy (\ref{cond1}) and
(\ref{cond2}). Then one has obviously
\begin{equation}
 \label{eq:cn-formula}
  c_{n} = \sum_{k=1}^{n} c_{n}(k).
\end{equation}
We have $c_{n}(k) = 0$ for $2^{k} \leq n$, because a polytope with at
most $n$ vertices cannot have dimension $n$. Furthermore, we have
clearly $c_{n}(n) = 1$ for all $n \in \mathbb{N}$. 

As mentioned already in Example \ref{sec:example_pre_hyp}, a
0/1-polytope that satisfies (\ref{cond1}), (\ref{cond2}), and
$\abs{V(P)} = 2^{k}$ has among its vertices exactly $2^{k}-1$
subgroups of index 2. Hence by ignoring the groups $U_{i} = V(P) \cap
H_{i}$ for $1 \leq i \leq n$, we get 
\begin{cor}
  For $k \leq n < 2^{k}$ one has
  \begin{equation*}
    c_{n+1}(k) = c_{n}(k-1) + c_{n}(k)(2^{k}-n-1). 
  \end{equation*}
\end{cor}
The first few values are given in the Table \ref{tab:number_of_polytopes}.
\begin{table}[htpb]
  \centering
  \begin{tabular}{|l||c|c|c|c|c|c|c|c|c|}
    \hline
    n $\backslash$ k & 1 & 2 &3 &4 &5 &6 &7 &8 & $c_{n}$ \\
    \hline
    \hline
    1 &  1 &   & & & & & & & 1 \\
    2 &  0 & 1 &   & & & & & & 1 \\
    3 &  0 & 1 & 1 &   & & & & & 2 \\
    4 &  0 & 0 & 5 & 1 &  & & & & 6 \\
    5 &  0 & 0 & 15& 16&  1 &  & & & 32 \\
    6 &  0 & 0 & 30& 175& 42& 1& & & 248 \\
    7 &  0 & 0 & 30& 1605& 1225& 99& 1 & & 2960\\
    8 &  0 & 0 & 0 & 12870&  31005& 6769& 219& 1& 50864 \\
    \hline
  \end{tabular}
  \caption{The number of $n$-dimensional 0/1 polytopes with $2^k$
    vertices that form a group. }
  \label{tab:number_of_polytopes}
\end{table}
It is easy to compute this number also for larger values of $n$. For
instance
\begin{gather*}
  c_{28} = 718897730072178204358180468879825453986397667929112558174208 \\
    c_{100} \approx 2.77 \cdot 10^{644}
\end{gather*}
Finally, using the Corollary we can show that, among the full
dimensional 0/1-polytopes with $2^{k}$ vertices the convex hulls of
linear codes are exceptional. For $1\leq k \leq n$, let $d_{n}(k)$
denote the number of all 0/1 polytopes with $2^{k}$ vertices
satisfying only condition (\ref{cond2}). Hence, the number $d_{n}$ of
all 0/1 polytopes of dimension $n$ trivially satisfies
\begin{equation}
  \label{eq:dn-formula}
  d_{n}\geq \sum_{k=1}^{n}d_{n}(k).
\end{equation}

Moreover, we get
\begin{prop}
  \begin{enumerate}[(i)]
  \item For $4\leq n < 2^{k} < 2^{n}$, one has
    \begin{equation}
      \label{eq:dn_inequality}
      d_{n}(k) \geq 2^{k}(2^{n} - 2^{k})c_{n}(k) > n2^{n-1}c_{n}(k).
    \end{equation}
  \item We have
    \begin{equation*}
      \lim_{n\to\infty} \frac{c_{n}}{d_{n}} = 0.
    \end{equation*}
  \end{enumerate}
\end{prop}
\begin{proof}
  \begin{enumerate}[(i)]
  \item Suppose that $U$ is a proper subgroup of $(W_{n},\oplus)$ with
    $\dim (\conv(U)) =n$ and $\abs{U} = 2^{k}$. 

    If $U'$ is another subgroup of $(W_{n},\oplus)$ with $\abs{U'} =
    \abs{U}$, then we have
    \begin{equation*}
      \abs{U\cap U'} \leq 2^{k-1}< 2^{n-1}
    \end{equation*}
    and, hence,
    \begin{equation*}
      \abs{ U \setminus U'} \geq 2^{k-1} > \frac{n}{2} \geq 2.
    \end{equation*}
    This means
    \begin{equation}
      \label{eq:atleastthree}
      \abs{U \setminus U'} \geq 3.
    \end{equation}
    There are $2^{k}(2^{n}-2^{k})$ subsets $V$ of $W_{n}$ with
    $\abs{V}=2^{k}$ and $\abs{V\setminus U} = \abs{U\setminus V} = 1$;
    namely, these are all sets of the form
    \begin{equation}
      \label{eq:setformincounting}
      V = (U\setminus\set{ u_{0}}) \cup \set{v_{0}} \quad\text{ with }
      u_{0}\in U, v_{0} \in W_{n}\setminus U.
    \end{equation}
    For $V$ as in \eqref{eq:setformincounting}, we get $\dim(\conv(V))
    = n$, because otherwise, $U\setminus\set{u_{0}}$ would be
    contained in a -unique- hyperplane $H$ with $v_{0}\in H$, a
    contradiction to $v_{0}\notin U$.  Together with
    \eqref{eq:atleastthree}, we obtain the first inequality in
    \eqref{eq:dn_inequality}. The second one is trivial in view of
    $2^{k}\leq 2^{n-1}$.
  \item By \eqref{eq:cn-formula}, \eqref{eq:dn-formula}, and
    \eqref{eq:dn_inequality} we get for $n\geq 4$:
    \begin{align*}
      \frac{c_{n}}{d_{n}} & \leq 2 \frac{c_{n} - 1}{d_{n} -1} \\
      & \leq 2
      \left(
        \sum_{k=1}^{n-1}c_{n}(k)
      \right)
      \left(
        \sum_{k=1}^{n-1}d_{n} (k)   \right)^{-1} \\
      & \leq 2 (n2^{n-1})^{-1} = \frac{2^{2-n}}{n}.
    \end{align*}
    This proves the second statement. \qedhere
  \end{enumerate}
\end{proof}

As a concluding remark, we study the question of constructing a
statistical model from a given polytope. Assume $P\subset \mathbb{R}^n$
satisfies (\ref{cond1}) and (\ref{cond2}), when does it come from a
hierarchical model? To begin with, observe that the number $m$ of
vertices of $P$ is a power of 2, since it must divide the number of
vertices of the cube $W_n$. We write $m = 2^N$. By Theorem
\ref{sec:algothm}, we know that $P$ can be constructed using the
algorithm. It is constructed from the $N$-cube by applying several
steps of the second type in Algorithm \ref{sec:iterative_algorithm}.
Therefore, every such polytope corresponds to a subset of columns
$\set{1},\ldots,\set{N},A_{N+1},\ldots,A_s$ in the matrix of
coordinates, or equivalently to a pre-hypergraph. However, this
pre-hypergraph is not unique.  If we are given only a list of
vertices, then there are several possibilities to choose a generator
matrix in standard form. As an example consider the polytope
$\mathcal{F}_{2,3}$:
\begin{center}
\begin{tabular}{|c||c|c|c|c|c|c|}
  \hline
  $x$ & $\set{1}$ & $\set{2}$ & $\set{3}$ & $\set{1,2}$ & $\set{1,3}$ &
  $\set{2,3} $ \\
  \hline\hline
 $(0,0,0)$ & 0 & 0 & $0$ & 0 & 0 & 0 \\
 \hline
 $(0,0,1)$ & 0 & 0 & $1$ & 0 & 1 & 1 \\
 \hline
 $(0,1,0)$ & 0 & 1 & $0$ & 1 & 0 & 1 \\
 \hline
 $(0,1,1)$ & 0 & 1 & $1$ & 1 & 1 & 0 \\
 \hline
 $(1,0,0)$ & 1 & 0 & $0$ & 1 & 1 & 0 \\
 \hline
 $(1,0,1)$ & 1 & 0 & $1$ & 1 & 0 & 1 \\
 \hline
 $(1,1,0)$ & 1 & 1 & $0$ & 0 & 1 & 1 \\
 \hline
 $(1,1,1)$ & 1 & 1 & $1$ & 0 & 0 & 0 \\
 \hline
 \hline
\end{tabular}
\end{center}
In fact every 3 by 3
submatrix which is, after permuting rows and columns, the identity matrix gives
a generator matrix in standard form.
Obviously, there are several such choices. Here, for instance we can choose the
canonical basis corresponding to the atoms:
\begin{center}
\begin{tabular}{|c||c|c|c|c|c|c|}
  \hline
  $x$ & $\set{1}$ & $\set{2}$ & $\set{3}$ & $\set{1,2}$ & $\set{1,3}$ &
  $\set{2,3} $ \\
  \hline\hline
 $(1,0,0)$ & 1 & 0 & $0$ & 1 & 1 & 0 \\
 \hline
 $(0,1,0)$ & 0 & 1 & $0$ & 1 & 0 & 1 \\
 \hline
 $(0,0,1)$ & 0 & 0 & $1$ & 0 & 1 & 1 \\
 \hline
 \hline
\end{tabular}
\end{center}

On the other hand, we can also reorder the columns and choose 
\begin{center}
\begin{tabular}{|c||c|c|c|c|c|c|}
  \hline
  $x$ & $\set{3}$ & $\set{12}$ & $\set{13}$ & $\set{1}$ & $\set{2}$ &
  $\set{2,3} $ \\
  \hline\hline
 $(1,1,1)$ & 1 & 0 & $0$ & 1 & 1 & 0 \\
 \hline
 $(0,1,0)$ & 0 & 1 & $0$ & 0 & 1 & 1 \\
 \hline
 $(1,1,0)$ & 0 & 0 & $1$ & 1 & 1 & 1 \\
 \hline
 \hline
\end{tabular}
\end{center}

When the generator matrix is chosen, one can apply the method given in
Section \ref{sec:CodingTheory} to construct a pre-hypergraph. For our
first generator matrix we get back the hypergraph we started with, for
the second choice we read of the pre-hypergraph
\[
\mathcal{A}' \defas \set{\set{1'},\set{2'},\set{3'},\set{1',3'},\set{1',2',3'},
\set{2',3'}},
\]
where we have introduced new units $\set{1',2',3'}$ corresponding to
the first three columns. This shows the ambiguity due to the choice of
a particular generator matrix if only the code is given.

\section*{Acknowledgment}
\small Thomas Kahle is supported by the Volkswagen Foundation, Nihat
Ay is supported by the Santa Fe Institute. The authors thank Bernd
Sturmfels for many valuable comments on the subject. We are grateful
to the anonymous referees for their reports which led to significant
improvement of the paper.


\end{document}